\documentclass[11pt,fleqn,twoside]{article}
\usepackage{amsfonts,latexsym}
\makeatletter
\newcommand{\prava}{\footnotesize\it
\begin{flushright}
\begin{minipage}{18cm}
Copyright \copyright 1998 by S.K. Paul and A.R. Chowdhury
\end{minipage}
\end{flushright}}

\newcommand{\name}[1]{\begin{flushleft}
                       \LARGE \bf #1
                       \end{flushleft}\vspace{-3mm}}

\newcommand{\Author}[1]{\begin{flushleft}
                       \it #1 \end{flushleft}}

\newcommand{\Adress}[1]{\begin{flushleft}
                       \it #1 \end{flushleft}}

\newcommand{\Date}[1]{\begin{flushleft}
                      \small  \it #1 \end{flushleft}}

\newcommand{\ehkol}{Author \ name}
\newcommand{\ohkol}{Article \ name}
\renewcommand{\@evenhead}{
\hspace*{-3pt}\raisebox{-15pt}[\headheight][0pt]{\vbox{\hbox to \textwidth
{\thepage \hfil \ehkol}\vskip4pt \hrule}}}
\renewcommand{\@oddhead}{
\hspace*{-3pt}\raisebox{-15pt}[\headheight][0pt]{\vbox{\hbox to \textwidth
{\ohkol \hfil \thepage}\vskip4pt\hrule}}}
\renewcommand{\@evenfoot}{}
\renewcommand{\@oddfoot}{}

\newcommand{\be}{\begin{equation}}
\newcommand{\ee}{\end{equation}}
\newcommand{\ba}{\hspace*{-5pt}\begin{array}}
\newcommand{\ea}{\end{array}}
\newcommand{\p}{\partial}
\newcommand{\ds}{\displaystyle}
\makeatother

\begin{document}
\setcounter{page}{349}

\thispagestyle{empty}

\renewcommand{\ehkol}{S.K. Paul and  A. Roy Chowdhury}
\renewcommand{\ohkol}{On the Analytical Approach to the $N$-Fold
B\"acklund Transformation}

\begin{flushleft}
\footnotesize \sf
Journal of Nonlinear Mathematical Physics \qquad 1998, V.5, N~4,
\pageref{paul-fp}--\pageref{paul-lp}. \hfill {\sc Letter}
\end{flushleft}

\vspace{-5mm}

\renewcommand{\footnoterule}{}
{\renewcommand{\thefootnote}{}
 \footnote{\prava}}

\name{On the Analytical Approach to \\ the {\itshape N}-Fold B\"acklund
Transformation \\ of Davey-Stewartson Equation}\label{paul-fp}

\Author{S.K. PAUL and  A. ROY CHOWDHURY}

\Adress{High Energy Physics Division, Department of Physics,\\
Jadavpur University, Calcutta 700032, India}

\Date{Received March 12, 1998; Accepted May 31, 1998}

\begin{abstract}
\noindent
$N$-fold B\"acklund transformation for the Davey-Stewartson equation
is constructed by using the analytic structure of the Lax
eigenfunction in the complex eigenvalue plane. Explicit formulae can
be obtained for a specif\/ied value of $N$. Lastly it is shown how
generalized soliton solutions are generated from the trivial ones.
\end{abstract}

\noindent
{\bf Introduction:} Inverse scattering transform holds a central
place in the analysis of nonlinear integrable system in either (1+1)- or
(2+1)-dimensions~[1]. On the other hand it has been found that
explicit soliton solutions can also be obtained by the use of
B\"acklund transformations in a much easier way~[2]. These B\"acklund
transformations are also useful in proving the superposition formulae
for these solutions. There have been many attempts to construct
explicit $N$-soliton solutions for nonlinear integrable system either
by B\"acklund transformations or the inverse scattering method. A
separate and elegant approach was developed by Zakharov {\it et al}~[3]
which relied on the pole structure of the Lax eigenfunction and use
of projection operators. In this letter we have used an
approach similar to that of Zakharov {\it et al} but have generated a
formulae for the $N$-fold BT of the nonlinear f\/ield variables
occurring in the (2+1)-dimensional Davey-Stewartson equation~[4].
Our approach is very similar to that of gauge transformation
repeatedly applied to any particular seed solution. Lastly we
demonstrate how non-trivial solutions are generated by starting with
known trivial ones.

\medskip

\noindent
{\bf Formulation:}
The Davey-Stewartson equation under consideration can be written as
\be \label{paul:1}
\ba{l}
\ds ir_t+r_{xx}-r_{yy} +r(A_2-A_1)=0, \\[1mm]
\ds iq_t+q_{yy}-q_{xx} +q(A_1-A_2)=0, \\[2mm]
\ds A_{1x}=-\frac 12 (q_y r +r_y q),\qquad
\ds A_{2y}=-\frac 12 (q_x r +r_x q).
\ea
\ee

\noindent
Equation (\ref{paul:1}) is known to be a result of the consistency of the
operators $T_1$ and $T_2$ written as $[T_1,T_2]\Psi=0$, where
\[
\ba{l}
\ds T_1 \Psi =\left\{ 2 \left( \ba{cc} \p_x & 0\\ 0& \p_y\ea
\right) +\left( \ba{cc} 0 &q\\ r & 0 \ea \right) \right\}\Psi
=0,\\[4mm]
\ds T_2 \Psi =\frac 12 \left\{ \left(i \p_t +\p^2_x +\p_y^2\right) +
\left(\ba{cc} 0& q_x\\ r_y & 0\ea \right) +A\right\}\Psi =
-\frac{K^2}{2}\Psi,
\ea
\]
and $\ds A=\left(\ba{cc} A_1 & 0\\ 0& A_2\ea \right)$,
$K=\mbox{const}$.

To proceed further we set
\[
A_1 =f_{1y}, \qquad A_2=f_{2x}
\]
and
\[
\Psi =\Phi \exp \left\{ i\left(\alpha +\lambda^{-2}\right) x -i
\left(\beta -\lambda^{-2}\right) y\right\}.
\]
Whence the Lax pair becomes,
\be \label{paul:4}
M \Phi =U \Phi, \qquad \Phi_t=V \Phi
\ee
with
\be \label{paul:5}
M =\left( \ba{cc} \p_x & 0 \\ 0 & \p_y \ea \right), \qquad
U =\left( \ba{cc} -i\left(\alpha +\lambda^{-2}\right) & -q/2\\[1mm]
-r/2 & i\left(\beta -\lambda^{-2}\right) \ea \right)
\ee
and
\be \label{paul:6}
V=\left( \ba{cc} i\Lambda +Q(\p_x, \p_y)+if_{1y} & iq_x\\[1mm]
ir_y & i\Lambda +Q(\p_x, \p_y)+if_{2y} \ea \right).
\ee
With
\[
\ba{l}
\ds \Lambda =K^2 -\left(\alpha+\lambda^{-2}\right)^2
-\left(\beta-\lambda^{-2} \right)^2,\\[2mm]
\ds Q(\p_x,\p_y) =i\left(\p_x^2+\p_y^2\right) -2 \left\{
\left(\alpha+\lambda^{-2} \right) \p_x
-\left(\beta-\lambda^{-2}\right)\p_y \right\},
\ea
\]
we can construct particular Jost solutions, corresponding to
$q=q_0=\mbox{const}$, $r=r_0=\mbox{const}$, $A_1=A_{10}=f_{1y}^0$ and
$A_2=A_{20}=f_{2x}^0$ with $f_{1y}^0=f_{2x}^0=\mbox{const}$. This
particular eigenvector $\Phi_0$ turns out to be
\[
\hat{\Phi}_0 =\left( \ba{cc}
\exp(\theta_1 x+\chi_1 y+\xi_1 t)  &
\exp(\theta_2 x+\chi_2 y+\xi_2 t) \\[1mm]
m_0\exp(\theta_1 x+\chi_1 y+\xi_1 t) &
n_0\exp(\theta_2 x+\chi_2 y+\xi_2 t) \ea \right)
\]
with
\[
\ba{l}
\ds \theta_1=-i\left(\alpha+\lambda^{-2}\right)+am_0, \qquad
\theta_2=-i\left(\alpha+\lambda^{-2}\right)+an_0, \\[1mm]
\ds \chi_1=b/m_0+i\left(\beta-\lambda^{-2}\right), \qquad
\chi_2=b/n_0+i\left(\beta-\lambda^{-2}\right),
\ea
\]
$\xi_1$, $\xi_2$ are arbitrary complex constants, $m_0$, $n_0$ are
arbitrary constants and $m_0\not= n_0$.

Note that $\det \Phi_0\not= 0$ so that $\Phi_0^{-1}$ exists.

Now suppose that $\Phi_{n-1}$ denotes the Lax eigenfunction,
corresponding to the $(n-1)$ soliton solution, and $B_n(x,y,t)$ be the
transformation which yields the $\Phi_n$ (solution corresponding to
the $n$ soliton case when applied to $\Phi_{n-1}$), that is
\[
\Phi_n(x,y,t,\lambda) =B_n(x,y,t,\lambda)\Phi_{n-1}(x,y,t,\lambda).
\]
Using the above Lax equations (\ref{paul:4}) and
(\ref{paul:6}), we can at once deduce the
equations satisf\/ied by $B_n$
\be \label{paul:9}
MB_n=U_n B_n-B_n U_{n-1},\qquad
\p_t B_n =V_n B_n -B_n V_{n-1},
\ee
where $U_n$, $V_n$ denote the Lax matrices.

Corresponding to the $n$-soliton solution: Note that $U$, $V$ are
even functions of $\lambda$, so that we can assume that
\[
B_n(-\lambda)=B_n(\lambda).
\]
We now assume $B_n$ to have simple pole structure in the complex
$\lambda$-plane, so that
\be \label{paul:11}
B_n(x,y,t,\lambda)=Q_n +\frac{2\lambda_n}{\lambda^2-\lambda_n^2}P_n.
\ee
We also assume that
\be\label{paul:12}
B^{-1}_n(x,y,t,\lambda)=Q_n' +\frac{2\lambda_n'}{\lambda^2-\lambda_n^2}P_n',
\ee
where
$P_n$, $Q_n$, $P_n'$, $Q_n'$ are matrix functions of $(x,y,t)$. The
condition $B_nB_n^{-1}= B_n^{-1} B_n=I$ leads to
\[
\ba{l}
B_n\left(\lambda_n'\right)P_n'=0, \qquad
P_nB_n^{-1}(\lambda_n)=0,\\[1mm]
B_n^{-1}(\lambda_n)P_n=0, \qquad P_n'B_n\left(\lambda_n'\right)=0.
\ea
\]

\noindent
{\bf Calculation of the matrices $Q_n$, $P_n$, $Q_n'$, $P_n'$:}
Let us now go back to equation (\ref{paul:9})
and use the expression (\ref{paul:11}) and
(\ref{paul:12}). Rewriting $U_n$ as
\[
U_n=-i\lambda^{-2} I+U_n'
\]
and using the form of $B_n$ given in (\ref{paul:11}) we get:
\[
\ba{l}
\ds MQ_n +\frac{2\lambda_n}{\lambda^2-\lambda_n^2} (MP_n) =
\left(-i\lambda^{-2}I+U_n'\right)
\left(Q_n +\frac{2\lambda_n}{\lambda^2-\lambda_n^2} P_n\right)\\[3mm]
\ds \qquad \qquad - \left(
Q_n+\frac{2\lambda_n}{\lambda^2-\lambda_n^2}P_n\right)
\left(-i\lambda^{-2}I+U_{n-1}' \right)
\ea
\]
which yields equations satisf\/ied by $P_n$ and $Q_n$:
\[
\ba{l}
M(P_n\Phi_{n-1}(\lambda_n)) =U_n(\lambda_n)
(P_n\Phi_{n-1}(\lambda_n)),\\[1mm]
M(Q_n\Phi_{n-1}(\lambda_n)) =U_n(\lambda_n)
(Q_n\Phi_{n-1}(\lambda_n)).
\ea
\]
Using the same form of $V_n$ as give in equation (\ref{paul:5})
in the time part, we get the following equations:
\[
\ba{l}
\ds \p_x\left( Q_n -\frac{2}{\lambda_n} P_n\right) =-\p_y\left(
Q_n-\frac{2}{\lambda_n} P_n\right), \\[3mm]
\ds\p_t(P_n\Phi_{n-1}(\lambda_n))
=V_n(\lambda_n)(P_n\Phi_{n-1}(\lambda_n)),\\[1mm]
\ds Q_{nt} =OP_2Q_n+2i(Q_{nx}\p_x+Q_{ny}\p_y)+D_nQ_n-Q_nD_{n-1},
\ea
\]
where operator $OP_2\equiv i\left(\p_x^2+\p_y^2\right)-
2(\alpha\p_x-\beta\p_y)$, $\alpha,\beta=\mbox{const.}$, and
$\ds D_n=\left(\ba{cc} if_{1y}^n & iq_{nx}\\
ir_{ny} & if_{2x}^n\ea \right)$.

As per the ansatz of Zakharov we search for $P_n$ in the form,
\[
P_n=\left(\begin{array}{c} \gamma_{n1} \\ \gamma_{n2} \end{array} \right)
(\delta_{n1}, \delta_{n2}).
\]
It is interesting to observe that
\[
(\delta_{n1}, \delta_{n2}) =(a_{n1}, a_{a2}) \Phi_{n-1}^{-1}(n),
\]
where $a_{n1}$, $a_{n2}$ are practically two constants.

Similarly for $P_n'$ and $Q_n'$, we set
\[
P_n'=\left(\begin{array}{c} \gamma_{n1}' \\[1mm] \gamma_{n2}'\end{array}
\right) \left(\delta_{n1}', \delta_{n2}'\right) \quad \mbox{and}
\quad Q_n'=\left( \ba{cc} \alpha_n' & \alpha_n''\\[1mm]
\beta_n' & \beta_n''\ea \right).
\]
Whence we get
\[
B_l(\lambda_l') =\left(
\ba{cc}
F_l^{11} +\sigma_l \gamma_{l1} \delta_{l1} &
F_l^{12}+\sigma_l\gamma_{l2}\delta_{l2}\\[1mm]
F_l^{21} +\sigma_l \gamma_{l2} \delta_{l1} &
F_l^{22}+\sigma_l\gamma_{l2}\delta_{l2}\ea \right),
\]
where
\[
\ba{l}
\ds F_l^{11}=f_l^{11}(t) \exp\{im_{1l}(x-y)\}, \qquad
F_l^{12}=f_l^{12}(t) \exp\{im_{1l}'(x-y)\}, \\[1mm]
\ds F_l^{21}=f_l^{21}(t) \exp\{im_{2l}(x-y)\}, \qquad
F_l^{22}=f_l^{22}(t) \exp\{im_{2l}'(x-y)\},
\ea
\]
and
\[
\ba{l}
\ds\sigma_l=\frac{2}{\lambda_l}-\frac{2\lambda_l}{\lambda_l^2-
\lambda_l^{\prime2}},\\[4mm]
\ds \delta_{l1}=-\frac{\delta_{l1}'F_l^{11}+\delta_{l2}'F_l^{21}}
{\sigma_l(\delta_{l1}'\gamma_{l1} +\delta_{l2}'\gamma_{l2})},
\qquad
\delta_{l2}=-\frac{\delta_{l1}'F_l^{12}+\delta_{l2}'F_l^{22}}
{\sigma_l(\delta_{l1}'\gamma_{l1} +\delta_{l2}'\gamma_{l2})},\\[4mm]
\ds \gamma_{l1}=-\frac{\gamma_{l1}'F_l^{11}+\gamma_{l2}'F_l^{12}}
{\sigma_l(\gamma_{l1}'\delta_{l1} +\gamma_{l2}'\delta_{l2})},
\qquad
\gamma_{l2}=-\frac{\gamma_{l1}'F_l^{21}+\gamma_{l2}'F_l^{22}}
{\sigma_l(\gamma_{l1}'\delta_{l1} +\gamma_{l2}'\delta_{l2})}.
\ea
\]
\[
B_l^{-1} (\lambda_l) =
\left[ \ba{cc} \alpha_l'+\varepsilon_l'\gamma_{l1}'\delta_{l1}' &
\alpha_l''+\varepsilon_l'\gamma_{l1}'\delta_{l2}' \\[1mm]
\beta_l'+\varepsilon_l'\gamma_{l2}'\delta_{l1}' &
\beta_l''+\varepsilon_l'\gamma_{l2}'\delta_{l2}'
\ea\right].
\]
Here$\ds
\varepsilon_l'=\frac{2\lambda_l'}{\lambda_l^2-\lambda_l^{\prime 2}}$,
along with
\[
\ba{l}
\ds
\delta_{l1}' =-\frac{\alpha_l'\delta_{l1} +\beta_l' \delta_{l2}}
{\varepsilon_l'(\delta_{l1} \gamma_{l1}' +\delta_{l2} \gamma_{l2}')},
\qquad
\delta_{l2}' =-\frac{\alpha_l''\delta_{l1} +\beta_l'' \delta_{l2}}
{\varepsilon_l'(\delta_{l1} \gamma_{l1}' +\delta_{l2}
\gamma_{l2}')},\\[4mm]
\ds
\gamma_{l1}' =-\frac{\gamma_{l1}\alpha_l' +\gamma_{l2}\alpha_l''}
{\varepsilon_l'(\gamma_{l1} \delta_{l1}' +\gamma_{l2} \delta_{l2}')},
\qquad
\gamma_{l2}' =-\frac{\gamma_{l1}\beta_l' +\gamma_{l2}\beta_l''}
{\varepsilon_l'(\gamma_{l1} \delta_{l1}' +\gamma_{l2}
\delta_{l2}')}.
\ea
\]
Finally the matrix $Q_l$ is given as
\[
Q_l=\left( \ba{cc}
\ds f_l^{11}(t) \exp\{im_{1l}(x-y)\}+\frac{2}{\lambda_l}\gamma_{l1}
\delta_{l1} &
\ds f_l^{12}(t) \exp\{im_{1l}'(x-y)\}+\frac{2}{\lambda_l}\gamma_{l1}
\delta_{l2}\\[4mm]
\ds f_l^{21}(t) \exp\{im_{2l}(x-y)\}+\frac{2}{\lambda_l}\gamma_{l2}
\delta_{l2} &
\ds f_l^{22}(t) \exp\{im_{2l}'(x-y)\}+\frac{2}{\lambda_l}\gamma_{l2}
\delta_{l2}
\ea \right)\!.
\]
It is also very convenient to rewrite the matrix elements of $Q$ and
$P$ in terms of Lax eigenfunctions $\Phi$. We collect these results
below without giving the detailed derivation,
\[
\ba{l}
\ds Q_l^{11} =\frac{R_l\Phi_{l-1}^{21} (\lambda_l)+
R_l'\Phi_{l-1}^{22}(\lambda_l)}{a_l\Phi_{l-1}^{22} (\lambda_l)-
b_l\Phi_{l-1}^{21} (\lambda_l)}, \qquad
Q_l^{12} =\frac{M_l\Phi_{l-1}^{21} (\lambda_l)+
M_l'\Phi_{l-1}^{22}(\lambda_l)}{a_l\Phi_{l-1}^{22} (\lambda_l)-
b_l\Phi_{l-1}^{21} (\lambda_l)},\\[5mm]
\ds
Q_l^{21} =\frac{L_l\Phi_{l-1}^{21} (\lambda_l)+
L_l'\Phi_{l-1}^{22}(\lambda_l)}{a_l\Phi_{l-1}^{22} (\lambda_l)-
b_l\Phi_{l-1}^{21} (\lambda_l)},
\qquad
Q_l^{22} =\frac{N_l\Phi_{l-1}^{21} (\lambda_l)+
N_l'\Phi_{l-1}^{22}(\lambda_l)}{a_l\Phi_{l-1}^{22} (\lambda_l)-
b_l\Phi_{l-1}^{21} (\lambda_l)},
\ea
\]
where
\[
\ba{l}
\ds R_l=-\frac{\lambda_l^{\prime 2} b_lf_l^{11}(t)\exp\{im_{1l}(x-y)\}}
{\lambda_l^{\prime2}},
\\[5mm]
\ds R_l'=\frac{1}{\lambda_l^{\prime 2}}
\left\{\lambda_l^2 a_lf_l^{11}(t)\exp\{im_{1l}(x-y)\} +
\left(\lambda_l^2-\lambda_l^{\prime 2}\right) b_l f_l^{12}(t)
\exp\{im_{1l}'(x-y)\}\right\},
\\[5mm]
\ds M_l=\frac{1}{\lambda_l^{\prime 2}}
\left\{-\left(\lambda_l^2 -\lambda_l^{\prime 2}\right)
a_lf_l^{11}(t)\exp\{im_{1l}(x-y)\} -\lambda_l^2 b_l f_l^{12}(t)
\exp\{im_{1l}'(x-y)\}\right\},
\\[5mm]
\ds M_l'=\frac{\lambda_l^{\prime 2} a_lf_l^{12}(t)\exp\{im_{1l}'(x-y)\}}
{\lambda_l^{\prime2}},
\qquad
\ds
L_l=-\frac{\lambda_l^{\prime 2} b_lf_l^{21}(t)\exp\{im_{2l}(x-y)\}}
{\lambda_l^{\prime2}},
\\[5mm]
\ds L_l'=\frac{1}{\lambda_l^{\prime 2}}
\left\{\lambda_l^2 a_l f_l^{21}(t)\exp\{im_{2l}(x-y)\} +
\left(\lambda_l^2-\lambda_l^{\prime 2}\right) b_l f_l^{22}(t)
\exp\{im_{2l}'(x-y)\}\right\},
\\[5mm]
\ds N_l=\frac{1}{\lambda_l^{\prime 2}}
\left\{-\left(\lambda_l^2 -\lambda_l^{\prime 2}\right)
a_lf_l^{21}(t)\exp\{im_{2l}(x-y)\} -\lambda_l^2 b_l f_l^{22}(t)
\exp\{im_{2l}'(x-y)\}\right\},
\\[5mm]
\ds N_l'=\frac{\lambda_l^{\prime 2} a_lf_l^{22}(t)\exp\{im_{2l}(x-y)\}}
{\lambda_l^{\prime2}}.
\ea
\]
The elements of the $P_l$ matrix are:
\[
\ba{l}
\ds P_l^{11} =-\frac{F_l}{\sigma_l} \frac{\Phi_{l-1}^{22}(\lambda_l)}
{a_l \Phi_{l-1}^{22} (\lambda_l) -b_l\Phi_{l-1}^{21}(\lambda_l)},
\qquad
 P_l^{12} =\frac{F_l}{\sigma_l} \frac{\Phi_{l-1}^{21}(\lambda_l)}
{a_l \Phi_{l-1}^{22} (\lambda_l) -
b_l\Phi_{l-1}^{21}(\lambda_l)},
\\[5mm]
\ds P_l^{21} =-\frac{F_l'}{\sigma_l} \frac{\Phi_{l-1}^{22}(\lambda_l)}
{a_l \Phi_{l-1}^{22} (\lambda_l) -b_l\Phi_{l-1}^{21}(\lambda_l)},
\qquad
 P_l^{22} =\frac{F_l'}{\sigma_l} \frac{\Phi_{l-1}^{21}(\lambda_l)}
{a_l \Phi_{l-1}^{22} (\lambda_l) -b_l\Phi_{l-1}^{21}(\lambda_l)},
\ea
\]
with
\[
\ba{l}
\ds
F_l=a_lf_l^{11}(t) \exp\{im_{1l}(x-y)\} +b_lf_l^{12}(t)
\exp\{im_{1l}'(x-y)\}, \\[2mm]
F_l'=a_lf_l^{21}(t) \exp\{im_{2l}(x-y)\} +b_lf_l^{22}(t)
\exp\{im_{2l}'(x-y)\}, \\[3mm]
\ds \sigma_l =\frac{2\lambda_l^{\prime
2}}{\lambda_l\left(\lambda_l^{\prime 2} -\lambda_l^2\right)}.
\ea
\]

\noindent
{\bf Construction of the nonlinear f\/ields:}
Once the form of the matrices $P_l$ and $Q_l$ are determined we can
construct the matrix $B_l$, so that the Lax eigenfunction
$\Phi_l(\lambda)$ for the next stage can be determined from that of
the previous one via,
\[
\Phi_l(\lambda) =B_l(\lambda) \Phi_{l-1}(\lambda).
\]
These expressions are very complicated, so we just quote one of them
to display their structure. For example,
\be\label{paul:28}
\ba{l}
\Phi_l^{11} (\lambda) =N_l/D_l, \\[1mm]
\ds D_l=a_l\Phi_{l-1}^{21} (\lambda_l)
-b_l\Phi_{l-1}^{22}(\lambda_l),\\[1mm]
\ds N_l=R_l\Phi_{l-1}^{21}(\lambda_l) \Phi_{l-1}^{11}(\lambda)+
\{R_l'-f_l(\lambda)F_l\} \Phi_{l-1}^{22}
(\lambda_l)\Phi_{l-1}^{11}(\lambda) \\[1mm]
\qquad
+\{M_l+f_l(\lambda)F_l\}\Phi_{l-1}^{21}(\lambda_l)\Phi_{l-1}^{21}(\lambda)
+M_l'\Phi_{l-1}^{22}(\lambda_l)\Phi_{l-1}^{21}(\lambda),
\ea
\ee
with similar expression for other elements $\Phi^{12}$, $\Phi^{21}$
and $\Phi^{22}$.

Now, for the determination of nonlinear f\/ields, consider
\[
U_n'Q_n =M_nQ_n+Q_nU_{n-1}',
\]
where
\[
M= \left( \ba{cc} \p_x & 0\\ 0& \p_y\ea \right), \qquad
U_n'=\left(\ba{cc} -i\alpha & \ds -\frac{q_n}{2}\\[2mm]
\ds -\frac{r_n}{2}& i\beta\ea \right), \qquad
U_{n-1}'=\left(\ba{cc} -i\alpha & \ds -\frac{q_{n-1}}{2}\\[2mm]
\ds -\frac{r_{n-1}}{2}& i\beta\ea \right).
\]
In the expression for $U_n'$ and $U_{n-1}'$ we take $\alpha=\beta=0$,
which at once yields
\be\label{paul:30}
q_n =-\frac{2Q_{nx}^{12}}{Q_n^{22}} +\frac{Q_n^{11}}{Q_n^{22}}q_{n-1}.
\ee
This is nothing but a simple recursion relation. Similarly,
\be \label{paul:31}
r_n =-\frac{2Q_{ny}^{21}}{Q_n^{11}} +\frac{Q_n^{22}}{Q_n^{11}}r_{n-1}.
\ee
Explicitly are can write,
\[
\ba{l}
\ds n=1, \quad q_1=-\frac{2Q_{1x}^{12}}{Q_1^{22}}
+\frac{Q_1^{11}}{Q_1^{22}}q_0,
\quad r_1=-\frac{2Q_{1y}^{21}}{Q_1^{11}}
+\frac{Q_1^{22}}{Q_1^{11}}r_0,
\\[4mm]
\ds n=2, \quad q_2 =-\frac{2Q_{2x}^{12}}{Q_2^{22}}-
\frac{2Q_2^{11} Q_{1x}^{12}}{Q_2^{22} Q_1^{22}}+
\frac{Q_2^{11}Q_1^{11}}{Q_2^{22} Q_1^{22}} q_0,\\[4mm]
\ds \phantom{n=2, \quad \ } r_2 =-\frac{2Q_{2y}^{21}}{Q_2^{11}}-
\frac{2Q_2^{22} Q_{1y}^{21}}{Q_2^{11} Q_1^{11}}+
\frac{Q_2^{22}Q_1^{22}}{Q_2^{11} Q_1^{11}} r_0.
\ea
\]
So far we have considered $f_l^{11}(t)$, $f_l^{12}(t)$,
$f_l^{21}(t)$, $f_l^{22}(t)$ to be functions of time or constants;
$m_{1l}$, $m_{1l}'$, $m_{2l}$, $m_{2l}'$ to be arbitrary constants;
$a_l$, $b_l$ to be arbitrary constants for all $l$. Now assume that
$f_l^{12}(t)= f_l^{21}(t)=0$ and $b_l=0$ for all $l$ values. So that
$R_l=M_l'=L_l=L_l'=N_l=F_l'=0$ for all $l$. In this case the form of
$B_l(\lambda)$ turns out to be:
$$
B_l(\lambda)=\!
\left(
\ba{cc} \frac{\lambda_l^2}{\lambda_l^{\prime 2}}\left(1-
\frac{\lambda_l^{\prime 2}-\lambda_l^2}{\lambda^2-\lambda_l^2}
\right) f_l^{11}(t)e^{im_{1l}(x-y)} &
\left(1+\frac{\lambda_l^2}{\lambda^2-\lambda_l^2}\right)
\left( 1-\frac{\lambda_l^2}{\lambda_l^{\prime 2}}\right) f_l^{11}(t)
e^{im_{1l}(x-y)} \theta_{l-1}\\[4mm]
\ds 0 &  f_l^{22}(t) e^{im_{2l}(x-y)}
\ea\!\! \right)\!\!,
$$
where $\ds
\theta_{l-1}=\frac{\Phi_{l-1}^{21}(\lambda_l)}{\Phi_{l-1}^{22}(\lambda_l)}$.

To write the formulae for the $n$-soliton solution in a compact form,
we note that
\[
\ds Q_n^{11} Q_{n-1}^{11} Q_{n-2}^{11}\ldots Q_R^{11}
 = A_{n-R+1}^n F_{n-R+1}^n(t) \exp\left\{
i\sum\limits_{j=1}^{n-R+1} m_{1,n-j+1} (x-y)\right\},
\]
where $A_{n-R+1}^n$, $F_{n-R+1}^n(t)$ stands for
\[
A_{n-R+1}= \prod_{l=1}^{n-R+1}
\frac{\lambda_{n-l+1}^2}{\lambda_{n-l+1}^{\prime 2}},
\qquad
F_{n-R+1}^n(t)=\prod_{l=1}^{n-R+1} f_{n-l+1}^{11}(t).
\]
Similar expressions can be written for $Q_n^{ij}$ and its products.
Using these, we at once obtain:
\[
\ba{l}
\ds q_n=\frac{-2\overline{m}_n Q_n^{12}}
{G_1^n(t) \exp\{im_{2,n}'(x-y)\}} +
\frac{A_n^n F_n^n(t) \exp\left\{ i\sum\limits_{j=1}^n
m_{1,n-j+1}(x-y)\right\} }{G_n^n(t)
\exp\left\{ i\sum\limits_{j=1}^n
m_{2,n-j+1}(x-y)\right\} }\\[5mm]
\ds \qquad -2\sum\limits_{K=1}^{n-1}
\frac{A_K^n F_K^n(t) \exp\left\{ i\sum\limits_{j=1}^K
m_{1(n-j+1)}(x-y)\right\} }
{G_{K+1}^n(t)
\exp\left\{ i\sum\limits_{j=1}^{K+1}
m_{2(n-j+1)}(x-y)\right\} }\overline{m}_{n-K}Q_{n-K}^{12},
\ea
\]
where we have set
\[
\overline{m}_0=a(m_0-n_0), \quad
G_K^n(t) =\prod_{l=1}^K f_{n-l+1}^{22} (t), \quad
\overline{m}_n' =-im_{1n}+\overline{m}_0',
\]
\[
\overline{m}_0' =b(1/m_0-1/n_0), \quad
F_k^n=\prod_{l=1}^k f_{n-l+1}^{11}(t), \quad
T_n=\frac{m_0}{n_0} \left( 1-\frac{\lambda_n^2}{\lambda_n^{\prime 2}}
\right),
\]
\[
Q_n^{12}=T_n f_n^{11}(t) \exp(im_{1n}x) \exp(\overline{m}_{0x}+
\overline{m}_n' y +\delta t).
\]
Here $\delta$ is a complex constant.


\medskip

\noindent
{\bf Discussions:}
In the above analysis we have demonstrated how the pole type ansatz
of Zakharov {\it et al} [3] can be used to generate a compact formulae for the
$N$-fold B\"acklund transformation in the case of the
Davey-Stewartson equation. The study yields two main results
exhibited in equations (\ref{paul:28}) and (\ref{paul:30}).
While the equation gives a
recursive procedure for the determination of the Lax eigenfunction
($\Phi_l$ corresponds to the $l$-soliton state) equation (\ref{paul:30}) and
(\ref{paul:31}) gives the corresponding recursion relation for the nonlinear
f\/ields. We have actually checked that for $n=1$ one obtains the one
soliton solution well known in the literature.

\label{paul-lp}

\end{document}